\documentstyle[12pt]{article}
\newtheorem{lem}{Lemma}[section]
\newtheorem{prop}{Proposition}

\newtheorem{defi}{Definition}[section]
\newtheorem{fact}{Fact}[section]
\newenvironment{proof}{{\bf Proof:}\newline}{\begin{flushright}$\Box$\end{flushright}}

\newcommand{\Cr}{{\bf Cr}}

\begin{document}

\title{Critical Circle Maps Near Bifurcation}

\author{Jacek Graczyk \\ Math Department, Warsaw University
\\ul. Banacha 2\\Warszawa 59, Poland  \and
Grzegorz \'{S}wia\c\negthinspace tek \thanks{ Partially 
supported by the NSF grant 431-3604A. 
This is a slightly modified version of the Stony Brook preprint with the same 
title. } \\Institute for Mathematical 
Sciences\\SUNY at Stony Brook\\Stony Brook, NY 11794, USA}
\date{July 15, 1991}
\maketitle
\begin{abstract}
We estimate harmonic scalings in the parameter space of 
a one-parameter family of critical circle maps. These estimates lead to
the conclusion that the Hausdorff dimension of the complement
of the frequency-locking set is less than $1$ but not less than $1/3$. 
Moreover, the rotation number is a H\"{o}lder continuous function of the 
parameter.

AMS subject code: 54H20   
\end{abstract}

\section{Preliminaries}
\subsection{Introduction}

This paper will present results about circle maps and  families of
circle maps that we were able to obtain during the past couple of years.
We will not discuss diffeomorphisms, which by far and large are the best
understood class of circle maps. In the present paper, we
will deal with critical homeomorphisms. Some methods and estimates can be
carried over to the non-invertible maps, but we only mention ~\cite{doktorat}
here.

We will provide rigorous proofs of our claims using only analytical tools and 
without resorting to the help of a computer (other than for typesetting this
text.) 

Let us start by defining the class of maps we consider.

\paragraph{The objects that we intend to investigate.}

Points of the real line can be projected onto the unit circle in the 
complex plane by means of the map
\[ x\rightarrow \exp(2\pi ix)    \; .\]
Maps from the real line project on the circle if they satisfy 
\[ F(x+1) - F(x) \in {\bf Z}\] for every real $x$. Obviously, for a 
continuous map this difference must be constant, and is the topological 
degree of the circle map. 

Unless necessary, we will not make a strong distinction between objects that
live on the circle and their lifts to the universal cover.   Whenever we want
to make a strong point of something being actually
on the circle, we will write (mod 1) near the formula.

If $x$ and $y$ are objects on the circle, $|x-y|$ is supposed to mean the 
distance in the natural metric.  
\subparagraph{Hypotheses.}
\begin{em}
We consider a family a circle maps given by
\[ F_{t}(x) = F_{0}(x) + t \; (mod\, 1)\]
where $t$ is a real parameter which ranges on the interval $(0,1)$ and 
$F_{0}\; (mod\, 1)$ is a degree one circle homeomorphism which in addition 
satisfies:
\begin{itemize}
\item
It is at least three times continuously differentiable.
\item
The derivative vanishes in exactly one point which is identified with $0$.
Also, $0$ is fixed by $F_{0}$.
\item
The function is differentiable enough times so as to satisfy 
\[\frac{dF_{0}^{l}}{dx^{l}}(0) \neq 0\]
for some $l$ where the $l$-th derivative exists and is continuous everywhere.
\end{itemize}  

In addition, we consider the corresponding family $f_{t}$ of lifts of maps
$F_{t}$. We denote with $\rho(t)$ the rotation number of $F_{t}$.
\footnote{We define the rotation number a little later.} 
\end{em}

\subparagraph{A generalization is possible.}
Most of our estimates would also work under a more general kind
of parameter dependence which only assumes that the family $f_{t}$ grows with
$t$ and there is some control over how the critical point moves. However, we
did not find that generalization worth the cost of  additional complications. 

\subparagraph{Theorem A.}
\begin{em}
Under the hypotheses listed above, the function $\rho(t)$ is H\"{o}lder 
continuous.
\end{em}

\subparagraph{Theorem B.}
\begin{em}
Consider the set
\[ \Omega' := \rho^{-1}({\bf R}\setminus {\bf Q}) \; .\]

Under our hypotheses the Hausdorff dimension of $\Omega'$ satisfies
\[ 1/3 \leq \mbox{HD}(\Omega') < 1\; .\]
\end{em}

\subparagraph{Remark B.}
\begin{em}
The upper estimate of Theorem B can be improved as follows:
there is a number $\alpha<1$, independent of the choice of $F_{0}$ which 
determines our family, so that 
\[ \mbox{HD}(\Omega') \leq \alpha\; .\]
\end{em}

While we do not formally prove Remark B, we give a justification for it.

\subparagraph{A comment on the results.}
Theorem A is the first result that we know about concerning regularity of
the rotation function for critical families. Theorem B is a refinement 
of the result that the measure of $\Omega'$ is zero . However, it still 
falls well short of numerically
established universality of the Hausdorff dimension (equal to about $0.87$
for cubic families). 

It worth noting that the lower estimate in Theorem B contradicts a certain 
conjecture based on extrapolating numerical data.  
The work \cite{phis} gave an 
asymptotic formula for the fractal dimension 
(which is essentially another name for Hausdorff dimension, see \cite{feder}) 
of $\Omega'$ which was expected
to  tend to $0$ as the critical exponent ($l$ in our notations) grew to 
infinity. This is, of course, contradicted by our result. In fact, the 
behavior of the family near bifurcation is what one tends to miss when doing 
numerics and one of our hopes is that our paper will help to place  
extrapolation of numerical results on a sounder basis.    

The results have three roots. One is bounded geometry of critical maps which
has been known for awhile. Another is much stronger estimates of the geometry 
near a bifurcation point which are new in the sharp form in which we use them.

Finally, there is a way of establishing similarity between objects in the
phase space and in parameters. We show an easy technique to achieve that
which, however, is much stronger than complicated and rather crude estimates
of this kind used in ~\cite{rat}.
    
\subsection{Topological Description of Dynamics}
The {\em rotation number} of a circle map $F$ is given by
\[ \rho(F) = \lim_{n\rightarrow\infty} \frac{f^{n}(x)}{n}\; (mod\, 1)\]
where $x$ is any point and $f$ any lift of $F$ to the universal cover. For 
the maps we consider, the limit always exists
and is independent of the choice of $x$ or a particular lift $f$. If the 
rotation number is 
irrational, it is a full topological invariant; and even if it is rational,
still a lot of information about the underlying dynamics can be read.

\paragraph{The structure on rotation numbers.}
There is no universal agreement on what is the best way to organize
rotation numbers. The reason for that may be that what ``best'' means may
depend on an author's particular objective. One way is to use the so-called
{\em Farey tree}, and another is based on {\em continued fractions}. We will 
base our approach on the Farey tree structure and only comment marginally 
on the connection with continued fractions. We will use lots of properties
of the rotation numbers in conjunction with the dynamics without proofs. 
All proofs would be elementary. Clues to the arithmetic part of the issue
can be found in ~\cite{hardy}, while the dynamical properties are scattered
throughout the literature. Some proofs and hints for others can be found in 
~\cite{rat}.

\paragraph{Farey trees.}
\begin{defi}\label{defi:n1,1}
We define the structure of a directed graph whose vertices are exactly
all rational numbers from $(0,1)$.

By definition, each vertex $p/q$ has exactly two outcoming edges. One leads to 
a smaller number called the ``left daughter'', and the other to a greater
number called the ``right daughter''.

If $p/q$ is in the lowest terms, we determine $u<p/q<v$ defined as the closest
neighbors in $[0,1]$ with denominators not larger than $q$. Then, the left 
daughter
is the rational number with the smallest denominator contained in $(u,p/q)$, 
while
the right daughter is the number with the smallest denominator contained in
$(p/q, v)$ (they are unique.) 

It is known that the graph of this relation is a connected binary tree and
$1/2$ is the root. This tree is called the {\em Farey tree}.
\end{defi}
\subparagraph{Coding.}
Thus, there is a one-to-one correspondence between a rational number mod 1
and a finite symbolic sequence $(a_{i})$ (the symbols are $L$ for the left 
daughter and
$R$ for the right daughter) which tells us how to go to the number from the 
root.

Also, there is a unique infinite symbolic representation of every irrational
number defined by the property that the rationals which correspond to finite 
initial segments of the code tend to the irrational number. The reader may try
to compute the value of the alternating code $LRLR\ldots$. 
\subparagraph{The degree.}
For every symbolic sequence, we define inductively its {\em turning points}
$m_{i}$. The first turning point $m_{1}$ is defined to be the least $i$ so
that $a_{i} \neq a_{i+1}$. If it does not exist, there are no turning points.
Once $m_{j}$ has been found, $m_{j+1}$ is the least $i>m_{j}+1$ so that
$a_{i} \neq a_{i+1}$. Again, if it does not exist, the turning point
sequence ends. This definition implies that if what one would intuitively call
a ``turning point'' occurs immediately after another turning point, it does
not count as a turning point. This is required for compatibility with the 
``harmonic'' description given a little later.

{\em The degree} of a rational number $u$ (denoted $\mbox{deg}(u)$) is the 
number of turning points of the corresponding symbolic sequence plus one.

The reader is kindly asked to determine the degree of the number coded by
\footnote{Five is the answer.}
\[ LRLLLRRLRLR\; . \]

\subparagraph{Closest returns.}
Consider an infinite symbolic sequence $\cal A$ which codes an irrational 
number $\rho$. We consider the sequence $q_{i}$ defined as denominators 
of consecutive rationals which correspond to symbolic sequences
\[ a_{1},\ldots, a_{t_{i}} \] 
where $t_{i}$ are turning points of ${\cal A}$. 

This sequence has a transparent interpretation in terms of the dynamics
of the rotation by $2\pi\rho$. Namely, fix a point on the circle and consider
the sequence of iterates which map this point closer to itself than any 
previous iterate. This turns out to be exactly the sequence $q_{i}$ defined 
above.
Because of this interpretation we will refer to $q_{i}$ as the sequence of 
{\em closest returns} for $\rho$.  

\subparagraph{Farey domains and the harmonic subdivision.}
\begin{defi}\label{defi:n1,2}
The interval $(P/Q,P'/Q')\subset [0,1]$ is called a {\em Farey domain}
if and only if either there is an edge between $P/Q$ and $P'/Q'$ in the Farey
tree or it is one the three: $[0,1]$, $[0,1/2]$, or $[1/2,1]$.

Also, we will say that $P/Q$ and
$P'/Q'$ are {\em Farey neighbors} if and only if $(P/Q, P'/Q')$ is a Farey
domain.
\end{defi}

\begin{fact}\label{fa:n6,1}
If $P/Q$ and $P'/Q'$ are Farey neighbors, then $|PQ'-P'Q|=1$ and 
$1/2 \leq Q/Q' \leq 2$.
\end{fact}

For every Farey domain we consider a sequence $u_{n}$ where $n$ ranges over
all integers. If $n$ is positive,
\[ u_{n} := \frac{(n+1)P+P'}{(n+1)Q+Q'}\; . \]
For $n$ non-positive
\[ u_{n} := \frac{P+(1-n)P'}{Q+(1-n)P'} \; .\]

It is an elementary check that $(u_{n},u_{n+1})$ are all Farey domains. We will
call the collection of Farey domains of this form the {\em harmonic subdivision
} of $(P/Q, P'/Q')$. The numbers $u_{n}$ themselves can be called endpoints 
of the harmonic subdivision. 

\subparagraph{Harmonic coding.}
Start with the unit interval and consider its subdivision of level $1$, that is
its harmonic subdivision. Next, harmonically subdivide every domain of the
subdivision of the previous level and so on. The Farey domains obtained
on the $k$-th level will be called {\em fundamental domains} of level $k$.
Now, consider the set of symbols
$S(n,n+1)$ and $E(n)$ where $n$ is an integer. If $S(n,n+1)$ occupies the
$k$-th slot of the code, it refers to the interval $(u_{n},u_{n+1})$ of the
harmonic subdivision of the fundamental domain defined by the preceding part
of the code. 
So, a fundamental domain of
level $k$ can be coded by a sequence of symbols of type $S(\cdot,\cdot)$ of
length $k$. Endpoints of the domains of level $k$ have a similar coding which
ends with an $E(\cdot)$ symbol. Finally, irrational numbers can be coded by
infinite sequences of symbols $S(\cdot,\cdot)$. This coding will be called
{\em harmonic coding} and sometimes is more useful than the Farey coding.
There is a close correspondence between the symbols in the harmonic code
of a number and the coefficients of its continued fraction expansion. In
particular, irrational numbers of constant type (i.e. with all continued 
fraction coefficients uniformly bounded) correspond to harmonic codes which
consist of symbols $S(n,n+1)$ with $|n|$ uniformly bounded. 

\subparagraph{Miscellaneous properties.}
All numbers in the same fundamental domain of level $k$ have the same closest
returns up to $q_{k-1}$.

\paragraph{ What  all this means for our family of maps.}
Since our assumptions guarantee that the rotation number is a non-decreasing
function of the parameter, all those objects from the realm of rotation 
numbers can be transported back to the parameter space. So we will talk 
of Farey domains and harmonic subdivisions in the parameter space as well.

There is a caveat, though. The rational numbers form an insignificant 
countable subset of the set of rotation numbers, but their preimage in the
parameter space, called the {\em frequency-locking} set, is a huge set of full
measure, as it was demonstrated in ~\cite{rat}. What happens is that 
the rotation number does not always grow with the parameter. It makes stops
on all rational numbers but never on irrationals. 

\subsection{Harmonic Scalings}
With the comments we have made so far, we hope to have explained the main 
purpose of the paper as stated in the abstract. 

The main technical lemma of the paper will concern the harmonic scalings
in the parameter space. To explain the notion we have to go back to our
construction of the harmonic subdivision. The elements
of the subdivision accumulate to the endpoints of the parent Farey
domain. Exactly how fast their sizes decrease is the question of harmonic 
scalings. 

It has long been known that the decrease is governed by a cubic law.
The earliest mention we found in the literature is ~\cite{kaneko}.
The first mathematically rigorous work which established the result was
~\cite{jon}. However, the estimates were non-uniform, i.e. it was proved
that the scalings are indeed asymptotically cubic near every frequency-
locking interval, but no estimate on how long one should wait to see the
asymptotics take over in each particular case.

\paragraph{The saddle-node phenomenon.} 
To fix the notations, let us concentrate on a  Farey domain $(P/Q, P'/Q')$.

\subparagraph{We can assume $Q<Q'$.}
Indeed so, because the map is symmetric with respect to the choice of an 
orientation. More precisely, instead of our family $f_{t}$ we could consider 
a family $\phi_{t}$ given by
\[ \phi_{t}(x) = -f_{-t}(-x)\]. 

It is easy to check that this operation means changing the direction in the 
parameter space and the orientation on the circle. The rotation number
of $\phi_{-t}$ is going to be equal to $1-{\bf rot}(f_{t})$. The new family
$\phi_{t}$ still satisfies our assumptions, but because the rotation numbers
have been flipped around one half, so has been the Farey tree.

\subparagraph{What happens near the lower extreme of $\rho^{-1}(P/Q, P'/Q')$.}

Directly below $\rho^{-1}(P/Q, P'/Q')$ there  is a frequency-
locking interval which belongs to $P/Q$. The most 
interesting point for us is the upper boundary of this frequency-locking.
This parameter value will be denoted by $t_{0}$. According to the general
theory (see ~\cite{jon}) $f_{t_{0}}$ is structurally unstable, even within
the family. It still has a periodic point of period $Q$, which must be 
neutral, but any increase of the parameter value will destroy it, while any
decrease will create two periodic points in its place.  

The graph of $f_{t_{0}}$ is tangent to the diagonal so that all non-periodic
orbits are attracted to the neutral orbits. When the parameter value increases,
the graph is pushed up and a funnel opens between the graph and the diagonal.

\subparagraph{How to measure the scalings.}

Let us take a closer look at the situation for a parameter value $t$ just a 
little above $t_{0}$ and on the interval between $0$ and the nearest critical 
point
of $f^{Q}_{t}$ on the right of $0$. This critical point must be a preimage of
$0$ and we denote it $f^{-q}_{t}(0)$ where $q=Q`-Q$ is the previous closest
return common to all maps from $\rho^{-1}(P/Q, P'/Q') $.

Since its rotation number is a little greater than $P/Q$, $f^{Q}_{t}$ moves 
the points a little to the right. Thus, the critical point, for example, moves
to the right every $Q$ iterates. Finally, it will leave the interval
$(0,f^{-q}_{t}(0))$ and the number of steps it takes tells us exactly
which domain of the harmonic subdivision we are in. So one way to determine 
the scalings could be to measure the interval in the parameter space between
$t_{k}$ for which the image of $0$ by $f^{Q}_{t_{k}}$ hits $f^{-q}_{t_{k}}(0)$
after $k$ steps and $t_{k+1}$ where the same requires $k+1$ steps.

This is  the idea we will follow.

\subparagraph{The scalings near the upper endpoint of $\rho^{-1}(P/Q, P'/Q')$
must follow the same rules.}

This is a rather trivial reduction. We can consider the Farey domain 
$((P+P')/(Q+Q'), P'/Q')$ and then flip the Farey tree as described above.
What we get is a Farey domain $(1-P'/Q', 1-\frac{P+P'}{Q+Q'})$ and now what 
used to be the scalings near the top of $\rho^{-1}(P/Q, P'/Q')$ now 
is equal to corresponding scalings at  the bottom of the new domain. 

So, will only consider the scalings near the lower extreme of the Farey domain,
but the results will automatically extend to the upper scalings as well. 

\subparagraph{The crucial role of the funnel.}
The key observation made by the authors of the earlier works is that the 
decisive factor in estimating the scalings is the time it takes the image of
$0$ to go through the funnel. There are two reasons for that. The first is that
the image of $0$ spends most of its time in the funnel; the other is that as 
we consider the scalings in a very close proximity of the end of 
$\rho^{-1}(P/Q, P'/Q')$, the corresponding changes of the parameter
are so tiny 
that they only bring about minute modifications to the orbit of $0$ in the 
region away from the funnel. The main factor which effects the orbit is the 
change in the funnel clearance.

To prove what has been said here and study the effect of the funnel clearance
on the orbit was the main achievement of both ~\cite{jon} and ~\cite{gra}.
The work ~\cite{jon} studied this effect for the critical maps, but at that
moment it was very hard technically to get uniform estimates, in particular
independent of the degree of the Farey domain $(P/Q, P'/Q')$  . 
In the meanwhile, ~\cite{gra} provided estimates which were uniform in this 
sense, but
only applied to families of diffeomorphisms. Now we are finally able to give
uniform estimates for the critical maps as well. 

\subsection{Notations}
\paragraph{Uniform constants.}
Letters $K$ with a subscript will be reserved for ``uniform constants.''
If we claim a statement which involves such constants we mean precisely
that for each occurrence of a constant $K_{\cdot}$ a positive value
can be inserted which will make the statement true. The choice is uniform
in the sense that once $f$ has been fixed, there is a choice of values
which makes the statement true in all cases covered.

We do not claim that the same value should always be inserted for the same 
constant. To avoid complete confusion, we only promise that in each particular
statement we will mean to substitute the same value for all occurrences of a 
constant. Thus, a statement $K_{1} < K_{1}$ is formally true in our 
convention, but not in accord with our promise.

Finally, if $a < b < c < d$, the {\em cross-ratio} of these points is defined 
by
\[ \Cr (a,b,c,d) := \frac{|b-a||d-c|}{|c-a||d-b|} \; .\]

\paragraph{The rate of change of $f^{Q}$ depending on the parameter.}

\begin{prop}\label{prop:1,1}
We define a set $S\subset T\times {\bf S}^{1}$ as follows. The $t$ component
ranges between the lower extreme of the frequency-locking $P/Q$ and the 
lower endpoint of the frequency-locking of $P'/Q'$. Once $t$ has been fixed,
$x$ must belong to $[0, F^{-q-Q}_{t}(0)]$. Then, if $(t_{1},x_{1})$ and
$(t_{2},x_{2})$ belong to $S$, the ratio
\[ \frac{\frac{\partial f^{Q}}{\partial\tau}(t_{1},x_{1})}
        {\frac{\partial f^{Q}}{\partial\tau}(t_{2},x_{2})} \]
is uniformly bounded and bounded away from $0$.  
\end{prop}

\subparagraph{Comment.} This is a new estimate. It was apparently unknown 
when \cite{rat} was written, as that paper instead uses a very complicated and
roundabout method in order to obtain inequalities in the parameter 
space. The main source of a tremendous technical progress that
has been made in the past four years is the ``real Koebe lemma'' of 
~\cite{guke} and its 
derivatives like the Distortion Lemma that we use here.

\paragraph{The role of negative Schwarzian derivative.}
We do not want to assume that our function $f$ has negative Schwarzian
derivative. However, there is a remarkable, though not hard, fact that high
iterates of our functions already have negative Schwarzian. We will use this 
fact in some of our future estimates, which will, therefore, be valid only
for a large enough number of iterates.  

The idea that high iterates become negative 
Schwarzian maps is certainly not new and has been known to people working
in the field. However, we are not aware of any proof in the literature. One 
reason for that may be that it is unclear how to formulate this result in 
reasonable
generality. Our lemma does not pretend to be general, but the reader will
see from the proof that an analogous argument will work in many other 
situations.

\subparagraph{Rescaling.}
There are few uniform estimates on higher order derivatives for high iterates.
However, it is often possible to get estimates if the map is properly rescaled.
On the formal level that means that we take an arc of the circle and an 
iterate of the function which maps a part of this arc into the arc. Next, we
conjugate it affinely, usually so that the length of the arc becomes one. 
We will refer to this operation most frequently as to ``changing 
the unit of length''. In most cases that is all we want to say, and we will
preserve original notations.   

\begin{prop}\label{prop:1,3}
Assume  that the degree of $P/Q$ is larger than a uniform $K_{1}$.
Consider a parameter value $t>t_{0}$ from the Farey domain domain 
$u(P/Q,P'/Q')$ and change the unit of length so that the length of  
$(0,f^{-q}_{t}(0))$ becomes $1$. 

Then, $Sf^{Q}_{t} \leq -K_{2}$ on 
\[ (F^{-Q}_{t}(0), F^{-q}_{t}(0)) \setminus \{ 0\} \; .\]  
\end{prop}

\subparagraph{Corollary to Proposition ~\ref{prop:1,3}.}
\begin{fact}\label{fa:n3,2} 
Consider a map $F_{t}$ whose rotation number belongs to a fundamental 
domain of degree $l$. Assume that for some $i\leq q_{l}$, the iterate 
$F^{i}_{t}$ is does not have critical points on an open interval $J$. Then, 
$F^{i}_{t}$ on $J$ can be as 
\[ F^{i}_{t} = F^{i_{1}}_{t} \circ N \circ F^{i_{2}}_{t}\]
where $i_{1}$ and $i_{2}$ are uniformly bounded while the Schwarzian of
$N$ is negative.
\end{fact}
\begin{proof}
Choose $k$ so that for
\[ (\frac{p_{k}}{q_{k}}, \frac{p_{k-1}+p_{k}}{q_{k-1}+q_{k}}) \] 
Proposition ~\ref{prop:1,3} holds. This can be done in a uniform fashion. 

To preserve familiar notations, denote $Q=q_{k+1}$ and $q=q_{k}$.
Then, we consider intervals $I_{1}=(F^{-Q}_{t}(0), F^{-q-Q}_{t}(0))$
and $I_{2}=(F^{-q-Q}_{t}(0), F^{-q}_{t}(0))$. On $I_{1}\cup I_{2}$, the first
return time of $F_{t}$ is either the $Q$ on $I_{1}$, or $q$ on $I_{2}$. 

We pick $i_{1}$ as the first moment when the image of $J$ hits $I_{1}\cup
I_{2}$. Then, $i_{2}$ is the largest so that
$F^{i-i_{2}}_{t}(J)$ contains one of the following points: $0$, an endpoint
of $I_{1}$, or an endpoint of $I_{2}$. It is clear that $i_{1}$ and $i_{2}$
are uniformly bounded as $k$ was uniformly bounded. Then, all other iterates
can be accounted for by composing the pieces of the first return map on 
$I_{1}\cup I_{2}$, whose Schwarzian is negative by Proposition ~\ref
{prop:1,3}.
\end{proof}

\subparagraph{The nonlinearity lemma.}
For a function $h$, we introduce a quantity $nh:=h''/h'$, also called 
{\em nonlinearity} of $h$.   

\begin{prop}\label{prop:n2,1}
Choose a parameter $t$ from the Farey domain
and change the unit of length as in Proposition ~\ref{prop:1,3}, i.e. so that
$F^{-q}_{t}(0)$ is in the unit distance from $0$. Then, we obtain the 
following estimate on the rescaled map:

\[ nF^{Q}_{t} < K_{1} \] on $(F^{Q}_{t}(0), F^{Q-q}_{t}(0))$.

If, in addition, we also assume that the claim of Proposition ~\ref{prop:1,3} 
holds, we get two more estimates
\begin{itemize}
\item
\[SF^{Q}_{t} > -K_{2} \; , K_{2} > 0\] on the same interval,
\item
If $|F^{Q}_{t}(x)-x| < K_{3}$, then 
\[ nF^{Q}_{t}(x) > K_{2} > 0\; .\] 
\end{itemize}
\end{prop}

\subsection{Prerequisites}
We will review basic facts on which the proof will be based.

\paragraph{Bounded geometry.}
\begin{fact}\label{fa:1,1}
Let $F$ be a map from our family with the rotation number in $(P/Q, P'/Q')$
so that $0<q=Q'-Q$. There is a uniform positive constant $K_{1}$ for 
which the following estimates hold: 
\begin{itemize}
\item
For any $x$, 
\[K^{-1} < \frac{|x-F^{Q}(x)|}{|x-F^{Q}(x)|} < K_{1}\; ,\]
and the same holds with $Q$ replaced by $Q'$.
\item
For $x=F^{i}(0)\: ,\;i=-q,\ldots,q$,
\[K^{-1} < \frac{|x-F^{Q}(x)|}{|x-F^{q}(x)|} < K\; .\]
and the same holds with $Q$ replaced by $Q'$.
\end{itemize}
\end{fact}

Unfortunately, this fact belongs to the ``folk wisdom'', and there is no clear
reference to the proof. If the orbit of the critical point is periodic, Fact 
~\ref{fa:1,1} was proved in ~\cite{rat}. Then, it was shown by M. Herman how to
generalize the argument so as to include maps with irrational dynamics as
well (see \cite{michel}.)

\paragraph{The Distortion Lemma.}
\begin{fact}\label{fa:1,2}
Let $f$ belong to our family. Consider $f^{-m}$ for some $m>0$ restricted 
to an interval $J=(a,b)$ and rescaled so that $|a-b|$ becomes the unit of 
length. If the following are satisfied:
\begin{itemize}
\item
\[ \frac{df^{m}}{dx} > 0\; ,\]
\item
intervals $J$, $f(J)$, \ldots, $f^{m-1}(J)$ are all disjoint,
\end{itemize}
then for $x\in J$, the uniform estimate
\[ |\frac{f''(x)}{f'(x)}| < K_{1}\max(1/|a-x|, 1/|b-x|)\]
holds. 
\end{fact}
\begin{proof}
This estimate follows from the Uniform Bounded Nonlinearity Lemma of 
\cite{poincare}.
\end{proof}

\section{Proofs of Propositions}
\subsection{Proof of Proposition \protect\ref{prop:1,1}}
The key tool that we will use more than once in this paper is  
an approximate representation of the parameter derivative in terms of 
the lengths of dynamically defined intervals in the phase space.
  
\begin{lem}\label{lem:2,2}
If $(t,x)\in S$, then   
\[ |\log(\frac{\partial f^{Q}_{t}}{\partial t}(x))-
\log(\sum_{i=1}^{Q}\frac{|F^{q+Q}(x) - F^{Q}(x)|}{|f^{i}(x)-f^{i+q}(x)|})| 
\leq K_{1}\; .\]
\end{lem}
\begin{proof}
The time derivative has the form 
\[ \sum_{i=1}^{Q} \frac{df^{Q-i}}{dx}(f^{i}(x)) \; .\]
By the Distortion Lemma, we can replace the derivatives by the ratios of 
nearby small intervals with a bounded  error and that is all that the lemma 
claims. That the Distortion Lemma can be used to give uniform estimates can
be checked with bounded geometry. 
\end{proof}

\begin{lem}\label{lem:5,1}
Let $t_{1}$ and $t_{2}$ belong to a Farey domain $(P/Q,P'/Q')$ union 
with two adjacent frequency-lockings, and $x$ be any point. Then
\[K_{1}^{-1} < \frac{|F^{q}_{t_{1}}(x) - x|}{|F^{q}_{t_{2}}(x) - x|} < K_{1}
\; . \]
\end{lem}
\begin{proof}
Let us assume that $t_{1} < t_{2}$. We fixed the configuration so that 
the $q$-th iterate moves points a little to the left.
We only need to prove that 
$|F^{q}_{t}(x) - x|$ does not grow to much as $t$ decreases from $t_{2}$ to 
$t_{1}$. Suppose that indeed it does grow a lot, and we will show a bound.
There is the point $F^{3q}_{t_{2}}(0)$ somewhere in the distance from $x$ 
which is comparable to the distance from $x$ to $F^{q}_{t_{2}}(x)$. For 
combinatorial reasons, $F^{q}_{t_{1}}(x)$ must be on the right of 
$F^{3q}_{t_{2}}(x)$. Indeed, suppose that $x$ needs $k$ iterates of 
$F^{q}_{t_{2}}$ to overcome the previous closest return. Then, only $k-2$ or
less iterates of $F^{q}_{t_{2}}$ would be needed which clashes with our 
assumption about rotation numbers.
\end{proof}

\begin{lem}\label{lem:5,4}
Let $t$ belong to $\rho^{-1}(P/Q, P'/Q')$. Choose $x$ and $y$ so that
\[ y\in [x,F^{Q}_{t}(x)] \; .\]
Then, the ratio
\[ \frac{|F^{Q}_{t}(x)-x|}{|F^{Q}_{t}(y)-y|} \]
is uniformly bounded and bounded away from $0$.
\end{lem}
\begin{proof}
This follows by the first statement of bounded geometry. We observe the 
ordering of points $F^{-Q}_{t}(x), F^{-Q}_{t}(y), x, y, F^{Q}_{t}(x),
F^{Q}_{t}(y)$ which can only be preserved if the claim of the Lemma holds.
\end{proof}

Finally, we get this estimate:
\begin{lem}\label{lem:5,2}       
Let $t_{1}$ and $t_{2}$ be contained in $\rho^{-1}[P/Q,P'/Q')$. Then, the 
time derivatives of
$f^{Q}_{t_{1}}(0)$ and $f^{Q}_{t_{2}}(0)$ are comparable within uniform 
constants.
\end{lem}
\begin{proof}
We first use Lemma ~\ref{lem:2,2} to convert the parameter derivatives to sums
of ratios of intervals. Next, we need to show that any the lengths of any two 
corresponding intervals in both sums are comparable within uniform constants. 

Thus, consider $(F^{i}_{t_{1}}(0), F^{i+q}_{t_{1}}(0))$ and
$(F^{i}_{t_{2}}(0), F^{i+q}_{t_{2}}(0))$. First, we can use 
Lemma ~\ref{lem:5,1} to compare the lengths of
$(F^{i}_{t_{1}}(0), F^{i+q}_{t}(0))$ and
$(F^{i}_{t_{2}}(0), F^{i+q}_{t}(0))$ instead where $t$ is conveniently chosen
in $\rho^{-1}(P/Q, P'/Q')$. But those intervals overlap, so 
Lemma ~\ref{lem:5,4} concludes the proof.
\end{proof}

\paragraph{The proof of Proposition ~\ref{prop:1,1}}
In view of Lemma ~\ref{lem:5,2}, the only thing that remains to be shown 
is the uniformity of change with respect to $x$ with $t$ fixed. By Lemma
~\ref{lem:2,2} this comes down to estimating the ratios of 
\[ \frac{|F^{i+q}_{t}(x)-F^{i}_{t}(x)|}{|F^{i+q}_{t}(0)-F^{i}_{t}(0)|}\; . \]
We need $t$ in the preimage of the a Farey domain bounded by $p/q$ to conclude
the argument by Lemma ~\ref{lem:5,4}. Since the left daughter of $p/q$ is at 
least $P/Q$, we can do this unless $t\in \rho^{-1}(P/Q)$. If so, we can repeat 
the argument of Lemma ~\ref{lem:5,2} of first replacing $t$ and using Lemma
~\ref{lem:5,1}, and then Lemma ~\ref{lem:5,4}.

\subsection{Proofs of propositions about the Schwarzian}
\paragraph{Basic lemmas.}
The Schwarzian derivative of a $C^{3}$ local diffeomorphism is given by:
\[ Sf := (f''/f')' - \frac{1}{2}(f''/f')^{2}\; .\]
 
There is this remarkable formula for the Schwarzian of a composition:
\[ S(f\circ g) = Sf\circ g\cdot (g')^{2} + Sg\]
which for iterates of $f$ becomes:
\[ Sf^{n} = \sum_{i=0}^{n-1} (Sf\circ f^{i})\cdot ((f^{i})')^{2} \; .\]

Let $F$ be one map from our family. Consider an iterate $F^{n}$ on an interval
$(a,b)$. For any point $x\in (a,b)$ we look at $S_{+}F^{n}$ defined to be
\[ S_{+}F^{n} := \sum_{i=0}^{n-1} \max(0, SF\circ F^{i}) (F^{i})^{,2} \;
.\] In this situation a have a lemma:
\begin{lem}\label{lem:3,1}
If the following conditions are satisfied:
\begin{itemize}
\item
there is a larger interval $(a',b')\supset (a,b)$ so that the derivative of 
$F^{n}$ does not vanish on $(a',b')$,
\item
all intervals $(a,b)$, \ldots, $F^{n-1}(a,b)$ are disjoint,
\item
$|(a,b)|$ is taken as the unit of length,
\end{itemize}

then $S_{+}F^{n}$ on $(a,b)$ is bounded by a constant which only depends on
$\Cr(a',a,b,b')$ times 
\[ \max\{\frac {|F^{i}(a,b)|}{|C|} : 0\leq i \leq n-1\} \]
where $|C|$ is the length of the whole circle.
\end{lem}
\begin{proof}
We will now use the length of the circle as our unit. Then, we want to use
the affine map $A$ which maps the unit interval to $(a,b)$ and consider
\[S_{+}(F^{n}\circ A)\; .\]

By the Distortion Lemma, each derivative in the formula defining $S_{+}$ is
comparable with the square of the length of the corresponding image of 
$(a,b)$, and the error is bounded proportional to $\Cr(a',a,b,b')$. On the 
other hand, the Schwarzian derivatives are bounded. Since the sum of lengths 
of all images is less than $1$, the lemma follows.
\end{proof}

\paragraph{Proof of Proposition ~\ref{prop:1,3}}
For simplicity, we will denote $F:=F_{t}$. We will use the length of the 
circle as a unit.

First, we will consider $F^{q}$ on the interval 
$(F^{-Q}_{t}(0), F^{-q}_{t}(0))$. $F$ clearly has negative
Schwarzian derivative where defined which is bounded away from $0$. Then, 
$G:=F^{q-1}$ on 
the image satisfies the assumptions of Lemma ~\ref{lem:3,1} (which follows 
from bounded geometry.) If $A_{1}$ and $A_{2}$ are affine and map the unit 
interval onto $(0, F^{-q}(0))$ and on its image by $F$ respectively, we know
that 
\[ S_{+} (G\circ A_{2})\] is small, thus we need to bound the derivative
of $A_{2}^{-1}\circ F\circ A_{1}$, which clearly is bounded by something close
to $l$ (the critical point order.) The quantity $S_{+}(G)$ can be bounded 
by Lemma ~\ref{lem:3,1}. Since the lengths of images of $(0, F^{-q}(0))$ are
exponentially small in terms of the degree of $P/Q$,  the composition 
$G\circ F\circ A_{1}$ has negative Schwarzian
bounded away from zero provided that the degree is sufficient.
The same argument works for the composition of
$F$ followed by $F^{Q-q-1}$ on $F^{q}(0, F^{-q}(0))$. The Proposition follows.

\paragraph{Proof of Proposition ~\ref{prop:n2,1}.}
The first statement follows immediately from the Koebe distortion lemma 
and bounded geometry. 

The key to the other two statements is the identity
\[ Dng = Sg + 1/2 (ng)^{2}\; .\]
Since by the first claim $nf^{Q}_{t}$ is bounded, $Sf^{Q}_{t}$ must satisfy
some lower bound at least for some points from 
$(f^{Q}_{t}(0), f^{Q-q}_{t}(0))$, or the nonlinearity would be strongly 
decreasing, and, therefore, could not be uniformly bounded. However, the
values of $Sf^{Q}_{t}$ are comparable within uniform multiplicative factors
for all point of that interval. That can best be seen from the chain expansion
\[ Sf^{Q}_{t}(x) = \sum_{i=1}^{Q} Sf^{i}_{t}(f^{i-1}_{t}(x))
(Df^{i-1}_{t}(x))^{2}\; .\]

By rescaling using 
affine maps of bounded slope, we can assume that the $x$ is at $0$,
while $f^{-q}_{t_{0}}(0)$ is at $1$. We then consider the class $\cal G$ of
diffeomorphisms defined on $[0,1)$ which satisfy these conditions:
\begin{enumerate}
\item
Their Schwarzian derivatives are negative and bounded away from $0$ by some
$-\beta$.
\item
For any $g\in {\cal G}$, $g(0)=|F^{Q}_{t}(x) - x|$ and 
$\lim_{y\rightarrow 1} g(y)=1$.   
\item
If $g(x) \geq x$.
\end{enumerate}

We will prove that $ng(0)$ is uniformly bounded away from $0$ in $\cal G$.
The class $\cal G$ can identified with the set of pairs. One element of every
pair is as function bounded from above by $-\beta$, and the other is a 
non-negative number. The function gives the Schwarzian derivative, and the 
number can be taken as $ng(0)$. By the classical Schwarz result, these 
together with the second condition will
determine a unique function which is in $\cal G$ if and only if the last 
condition is satisfied. We then see that if $g_{1}, n$ is in $\cal G$, then 
$g_{2}, n$ is also in $\cal G$ provided that $g_{1} \leq g_{2} \leq -\beta$. 
That is because the corresponding solutions to the Schwarz problem will 
also satisfy the inequality. So, we can restrict our attention to maps
$-\beta, n$. By a compactness argument, unless  $n$ is uniformly bounded
away from $0$, there is a map in $\cal G$ with $n=0$. But it can be 
explicitly shown, that such a map cannot satisfy the third condition if
$K_{3}$ was chosen small enough in the statement of the Proposition.

\section{Scalings rules}
\subsection{Notations and the main result}
Throughout this chapter, we fix a Farey domain $(P/Q, P'/Q')$. As explained in
the Introduction, we may adopt the convention $0<q=Q'-Q$. 
We will denote
the frequency-locking interval which belongs to a rational number $u$ with 
$\alpha(u)$. The frequency-locking intervals which bound 
$\rho^{-1}(P/Q, P'/Q')$ from below
and from above are $\alpha(P/Q)$ and $\alpha(P'/Q')$. For every rational
 $u$, within $\alpha(u)$ there is a unique point $c(u)$, called the 
{\em center} of the corresponding frequency-locking, and characterized by the
property that the critical point is periodic. 

Next, we consider the sequence $u_{n}$ of endpoints of the harmonic 
subdivision of $(P/Q, P'/Q')$.

We define 
\[ J_{n} := (u_{n+1}, u_{n})\; ,\;\;\; \mbox{and} \]
\[ J := (c(P/Q), c(P'/Q')) \; .\] 

\begin{defi}\label{defi:n3,1}
{\em Harmonic scalings}, are the ratios
\[ h_{n} = \frac{|J_{n}|}{|J|} \]
for $n\in {\bf Z}$.
\end{defi}

Our main result about harmonic scalings is contained in the following 
proposition:
\begin{prop}\label{prop:3,1}
Harmonic scalings decrease no faster than according to a uniform cubic law, 
i.e. 
\[ h_{n} \geq \frac{K_{1}}{ |n|^{3}+1}\; , \]
where $K_{1}$ does not depend on $\cal D$.

If, in addition,  the claim of Proposition ~\ref{prop:1,3} holds, then  
\[h_{n} \leq \frac{K_{2}}{|n|^{3}+1}\;, .\] 
\end{prop}
  
The rest of this section will be devoted to the proof of Proposition 
~\ref{prop:3,1}.
\subsection{First estimates}
As it was noticed in section 1, it is enough to estimate 
$h_{n}$ for $n$ positive, since we can use the Farey domain 
$((P+P')/(Q+Q'), P'/Q')$ and then flip the Farey tree. We define 
\[ t_{n} := c(u_{n}) \]
and $t_{\infty} = c(P/Q)\; , \;\; t_{-\infty}=c(P'/Q')$.

As $t$ moves from $t_{n+1}$ to $t_{n}$, $F^{(n+1)Q}_{t}(0)$ travels from 
$F_{t_{n+1}}^{-q-Q}(0)$ to $F_{t_{n}}^{-q}(0)$. Thus,
\[ |J_{n}| \frac{df^{(n+1)Q}_{t|t=z}(0)}{dt} = 
|F_{t_{n+1}}^{-q-Q}(0) - F_{t_{n}}^{-q}(0)|\; \]
where $z$ is given by the Mean Value Theorem, thus $z\in (t_{n+1}, t_{n})$.
An analogous argument shows that 
\[ |J| \frac{df^{Q}_{t|t=w}(0)}{dt} = |F^{-q}_{t_{-\infty}}(0)-0| \; \]
this time with $w\in \alpha(P/Q)\cup \rho^{-1}(P/Q, P'/Q')$.   

These equalities allow us to express the harmonic scalings in terms of ratios
of time derivatives and lengths of relevant intervals. We will now work to make
this relation as simple as possible.

\begin{lem}\label{lem:n1,2}
The ratio of lengths 
\[\frac{|F_{t_{n+1}}^{-q-Q}(0) - F_{t_{n}}^{-q}(0)|}
{|F^{-q}_{t_{-\infty}}(0)-0|}\]
is uniformly bounded and bounded away from $0$.
\end{lem}
\begin{proof}
The tools of the simple proof are Lemma ~\ref{lem:5,1} and bounded geometry.
Since $t_{n+1}$ and $t_{n}$ are in frequency-lockings adjacent to a 
Farey domain $(u_{n+1}, u_{n})$ which belongs to the harmonic subdivision 
of $(P/Q, P'/Q')$, Lemma ~\ref{lem:5,1} can be applied with $Q$ to see that
$|F_{t_{n+1}}^{-q-Q}(0) - F_{t_{n}}^{-q}(0)|$ is uniformly comparable to
$|F_{t_{n}}^{-q-Q}(0) - F_{t_{n}}^{-q}(0)|$. By bounded geometry, that is in
a uniform relation to $|F_{t_{n}}^{-q}(0)-0|$, and another application of 
Lemma ~\ref{lem:5,1} with $q$ gives the claim.
\end{proof}

In view of Lemma ~\ref{lem:n1,2}, $h_{n}$ is uniformly comparable to
\begin{equation}\label{equ:n1,1}
 \frac{\frac{df^{Q}_{t|t=w}(0)}{dt}}{\frac{df^{(n+1)Q}_{t|t=z}(0)}{dt}}\; .
\end{equation}

We now concentrate on estimating expression ~\ref{equ:n1,1}.
\paragraph{Estimates of time derivatives.}

We first consider the denominator. By bounded geometry
and the Distortion Lemma, we see that  
\[ \frac{df^{(n+1)Q}_{t|t=z}(0)}{dt} \]
is uniformly comparable to
\[ \sum_{k=0}^{n-1}\frac{\partial f^{Q}}{\partial t}(z ,f^{kQ}_{z}(0)) 
\frac{|F^{nQ}_{z}(0)-F^{(n+1)Q}_{z}(0)|}
{|F^{kQ}_{z}(0)-F^{(k+1)Q}_{z}(0)|} \; .\]

Here, we replaced the spatial derivatives with ratios of intervals like in
the proof of Lemma ~\ref{lem:2,2}. Next, we use Proposition
~\ref{prop:1,1} to replace all parameter derivatives in this expression by 

\[ \frac{\partial f^{Q}}{\partial t}(z,0)\; , \]
again preserving comparability by uniform constants.

We see that expression ~\ref{equ:n1,1} is uniformly comparable to 
\[ \frac{ \frac{\partial f^{Q}}{\partial t}(w,0)}
{ \frac{\partial f^{Q}}{\partial t}(z,0)} \]

divided by 
\[\sum_{k=0}^{n-1}  \frac{|F^{nQ}_{z}(0)-F^{(n+1)Q}_{z}(0)|}
{|F^{kQ}_{z}(0)-F^{(k+1)Q}_{z}(0)|} \; \]

The first factor is uniformly bounded and bounded away from $0$ as a 
consequence of Lemma ~\ref{lem:5,2}.

The results of this section are summarized by the following proposition:

\begin{prop}\label{prop:n1,1}
For $n$ positive, the harmonic scaling $h_{n}$ is uniformly comparable to
\[\frac{1}{\sum_{k=0}^{n-1}  \frac{|F^{-q}_{z}(0)-0|}
{|F^{kQ}_{z}(0)-F^{(k+1)Q}_{z}(0)|}} \; .\]
\end{prop}
\begin{proof}
The proposition has almost been proven
apart from the fact that 
\[|F^{nQ}_{z}(0)-F^{(n+1)Q}_{z}(0)|\]
 is comparable
to $|F^{-q}_{z}(0)-0|$. That, however, follows from bounded geometry and 
Lemma ~\ref{lem:5,4} 
\end{proof}

\subsection{Saddle-node estimates}
In this section we will get uniform estimates for harmonic scalings depending
on $n$. To this end, we will use Proposition ~\ref{prop:n1,1} which makes this
task equivalent to estimating the geometrically given sum.

\paragraph{It is sufficient to consider $h_{n}$ with $n$ large.}
This is a simple corollary to proposition ~\ref{prop:n1,1}. The formula 
given there gives values bounded away from $0$ and infinity in a uniform
fashion for any bounded $n$. Thus, the claim of Proposition ~\ref{prop:3,1}
can be satisfied by choosing the uniform constants appropriately.  
\paragraph{Normalization.}

Consider $F^{Q}_{z}$ with $z\in J_{n}$ as in Proposition ~\ref{prop:n1,1}.
Consider $\chi\in (0, F^{-q}_{z}(0))$, which does not have to unique, so that 
\[ |F^{Q}_{z}(\chi) - \chi| = \inf \{ |F^{Q}_{z}(x) - x|\; : x\in 
(0, F^{-q}_{z}(0)) \} \; .\]

We change the coordinates by an affine map so that $\chi$ goes to $0$ and
$F^{-q}_{z}(0)$ goes to $1$. The critical point then lands at some point
$c$ whose distance from $0$ is uniformly bounded and bounded away from $0$. 
In these coordinates $F^{Q}_{z}$ becomes a map $\phi$. By the Distortion Lemma,
the second derivative of $\phi$ is bounded on a uniform neighborhood of $0$.

\paragraph{Approximation rules.}
We say that $\phi$ satisfies the $(\alpha,\kappa)$ {\em upper approximation rule}
if $\phi$ is not greater than the map 
\[ x\rightarrow x+\alpha x^{2} +\phi(0)\]
on some interval $(-\kappa,\kappa)$.

Analogously, $\phi$ satisfies the $(\alpha,\kappa)$ {\em lower approximation rule}
if there is the converse inequality.

Since the second derivative of $\phi$ is bounded on a uniform neighborhood
of $0$, there is a uniform choice of $(\alpha,\kappa)$ so that the upper 
approximation rule is always satisfied. In this sense we will say that $\phi$
satisfies the uniform upper approximation rule.

If the claim of Proposition ~\ref{prop:1,3} holds, there is also a uniform 
way to satisfy the lower approximation rule as a result of the estimate of
Proposition ~\ref{prop:n2,1}. So, we will consider maps which satisfy the
uniform lower approximation rule.

The approximating quadratic maps chosen in this uniform way will be called
$\phi_{u}$ for the map given by the upper approximation rule, and $\phi_{l}$
for the other map. 

The advantage of approximation rules is the orbits under quadratic maps can be
examined more or less explicitly and the interesting quantities simply
calculated. This was basically the idea of the authors who previously 
contributed to the subject, (see \cite{kaneko}, also \cite{jon} and 
\cite{gra}, which by no means 
exhaust the list as this trick was discovered independently a couple of 
times.) We will give only essential arguments and the reader will 
be able to find complete technical explanations in ~\cite{gra}. 

\subparagraph{Lemmas about quadratic maps.}

Consider a function $\Phi$
\[ x\rightarrow x+\alpha x^{2}+\epsilon\]
defined on     $(-\kappa, \kappa)$. Also assume that $\alpha > \beta > 0$ for 
some $\beta$, while $\kappa > \gamma > 0$ for some $\gamma$,
 and $\epsilon > 0$. A {\em maximal orbit} is a sequence  
$(y_{i})_{0\leq i \leq l}$
for which $y_{i+1}=\Phi(y_{i})$ and $y_{0}$ has no preimage while $y_{l}$ is no
longer in the domain. We have two facts about maximal sequences which are 
proved, though not explicitly stated,  in \cite{gra}. 

\begin{fact}\label{fa:n1,2}
The length of a maximal sequence $l$ and $\epsilon$ are related by
\[ K_{1}(\beta, \gamma)^{-1}\epsilon \leq \frac{\alpha}{l^{2}} 
\leq K_{1}(\beta,\gamma)\epsilon \]
where $K_{1}(\cdot,\cdot)$ is a positive function of $\beta,\gamma$ only. 
\end{fact}

\begin{fact}\label{fa:n1,3}
Let $l_{0}$ be the number of points $y_{i}$ from the maximal orbit which
satisfy $y_{i+1}-y_{i} < 2\epsilon$. Then
\[ \frac{l_{0}}{l} > K_{2}(\beta) \]
where again $K_{2}(\cdot)$ is a positive function of $\beta$ only.
\end{fact}

\subparagraph{Relating $\phi(0)$ and $n$.}    

\begin{lem}\label{lem:n2,2}
If $\phi$ satisfies the uniform upper approximation rule, 
\[ \frac{K_{1}^{-1}}{n^{2}} \leq \phi(0)\; . \]
\end{lem}
\begin{proof}
The number $n$ is at least the length of a maximal orbit by $\phi_{u}$.
The claim then follows by Fact ~\ref{fa:n1,2}.
\end{proof}     

\paragraph{Proof of the first estimate of Proposition ~\ref{prop:3,1}.}
There are $n$ terms in the denominator of the formula given by 
Proposition ~\ref{prop:n1,1}. By Lemma ~\ref{lem:n2,2}, they are at most of
the order of $n^{-2}$. The estimate on $h_{n}$ from below follows.

\paragraph{The negative Schwarzian case.}
The other estimate of Proposition ~\ref{prop:3,1} is harder, but since we
only claim it in the negative Schwarzian case, we can be aided by the strong
claim of Proposition ~\ref{prop:n2,1}. Among other things, we now know 
that $\phi$ satisfies both uniform approximation rules. Also, $\phi'$ has 
exactly one minimum which must be attained on the right of $0$. Consider 
$(-\kappa,\kappa)$ so that both uniform approximation rules hold.

\paragraph{Only what's inside $(-\kappa, \kappa)$ counts.}

We prove this lemma:
\begin{lem}\label{lem:n2,1} 
The number of iterates $k$ so that $0<k<n$ and 
$\phi^{k}(c)\notin (-\kappa, \kappa)$ is uniformly bounded. Moreover, for all such
values of $k$, $\phi^{k+1}(c) - \phi^{k}(c)$ is uniformly bounded away from 
$0$. 
\end{lem}
\begin{proof} 
Actually, the first part of the claim obviously follows from the second. 
This is only a problem if there are many iterates $k$ so that $\phi^{k}(c)
\in (-\kappa, \kappa)$, otherwise $\phi(0)$ is bounded away from $0$.   
Thus, by the approximation rules, if $k_{1}$ and $k_{2}$ are the smallest and 
the largest $k$ so that $\phi^{k}(c)\in (-\kappa, \kappa)$, the distances
$\phi_{k_{1}+1}(c) - \phi_{k_{1}}(c)$ and $\phi_{k_{2}+1}(c) -
\phi_{k_{2}}(c)$ are uniformly large. But they are still smaller than the 
analogous distances for the values of $k$ for which $\phi^{k}(c)\notin (-\kappa, \kappa)$
by the negative Schwarzian property.
\end{proof}

That means that in the sum of Proposition ~\ref{prop:n1,1} the contribution
from the terms corresponding to  values of $k$ with the property that 
$\phi^{k}(c)\notin (-\kappa, \kappa)$ is uniformly bounded. On the other hand,
it is 
clear that the whole expression grows at least as $n^{2}$. Thus, for values 
of $n$ uniformly sufficiently large, only the points of the orbit which are 
contained in $(-\kappa, \kappa)$ can be considered, and the result will 
approximate the 
whole sum up to a uniform multiplicative factor.

\paragraph{Essential estimates.}
Consider the smallest $k$ so that $\phi^{k}(c)\in (-\kappa, \kappa)$, and 
denote this 
point with $x_{0}$. Correspondingly, let $x_{l}$ be the last point of the
orbit still in $(-\kappa, \kappa)$, and in between we get a sequence which 
satisfies
$x_{i+1} = \phi(x_{i})$ for $i=0, \ldots, l-1$. We are interested in 
estimating the sum
\begin{equation}\label{equ:n1,2}
\sum_{i=0}^{l-1} \frac{1}{x_{i+1}-x_{i}} \; .
\end{equation}
from below.

\begin{lem}\label{lem:n4,2}
In our situation,
\[    \frac{K_{1}^{-1}}{n^{2}} \leq \phi(0) \leq \frac{K_{1}}{n^{2}}\; . \] 
\end{lem}

\begin{proof}
This is a stronger version of Lemma ~\ref{lem:n2,2} under stronger 
assumptions, and the proof is really the same.
\end{proof}

\subparagraph{Estimating the sum ~\ref{equ:n1,2} from below.}
Let $y$ be the largest point of the orbit of $c$ by $\phi$ still negative.
Then $(y^{l}_{i})$ and $(y^{u}_{i})$ denote maximal orbits for $\phi_{l}$ and
$\phi_{u}$ respectively which also contain $y$. The sum given by \ref{equ:n1,2}
is larger than the corresponding sum for $(y^{u}_{i})$, as the intervals 
occurring in the latter are fewer and longer. An analogous argument shows
that the sum for $(y^{l}_{i})$ bounds the interesting expression from above.

To estimate
\[ \sum_{i} \frac{1}{y^{u}_{i+1}-y^{u}_{i}} \]
from below we have to use Lemma ~\ref{lem:n4,2} which then asserts
simply that $\phi(0)$ is comparable to $n^{-2}$. Then, Fact ~\ref{fa:n1,2}
implies that $n$ is comparable to the length of the maximal orbit of both
$\phi_{u}$ and $\phi_{l}$.

Then, by Fact ~\ref{fa:n1,3} we see that the number of intervals for $\phi_{u}$
with lengths not greater than $2\phi(0)$ is still comparable with the
length of the maximal orbit, thus with $n$.

Thus, we  get a uniform cubic estimate from below. Again, Proposition 
~\ref{prop:n1,1} immediately enables us to derive the second claim of
Proposition ~\ref{prop:3,1}.

\section{H\"{o}lder continuity of the rotation number}

In this section, we will prove Theorem A announced in the Introduction.  
\subsection{Some consequences of scalings rules}
We will draw certain ``H\"{o}lder type'' estimates as consequences of 
Proposition ~\ref{prop:3,1}. We want to emphasize that we need estimates 
everywhere on the parameter space and we are unwilling to assume that the
denominators of Farey domains that we work with are large enough. So, only the 
first claim of Proposition ~\ref{prop:3,1} holds, i.e.
\[ h_{n} \geq \frac{K_{1,P.~\ref{prop:3,1}}}{|n|^{3}} \; .\]

\paragraph{Two estimates.}
We fix our attention on a Farey domain $(P/Q, P'/Q')$ subject to our usual 
convention $0 < q = Q' - Q$. We consider the harmonic subdivision by points
$u_{n}$. Also, the centers of mode-locking intervals are denoted with $t_{n}$
($n$ may be infinite) as in the Scalings Rules section.  
\begin{lem}\label{lem:n3,1}
Let $u=u_{n}\: ,\, v=u_{n+1}$. Then, 

\[|u-w||c(u)-c(w)|^{-\alpha}
\leq \beta|P/Q-P'/Q'||c(P/Q)-c(P'/Q')|^{-\alpha}\]
with uniform $0 <\alpha, \beta < 1$ 
\end{lem}
\begin{proof}
We have the following estimate 
as a consequence of Proposition ~\ref{prop:3,1}:

\[ \frac{|\rho(t_{n+1}) - \rho(t_{n})|}{|\rho(t_{-\infty})-\rho_{\infty}|}
\frac{|t_{-\infty}-t_{\infty}|^{\alpha}}{|t_{n+1}-t_{n}|^{\alpha}} \leq
\frac{n^{3\alpha} QQ'}{K_{1,P.\ref{prop:3,1}}^{\alpha} (nQ+Q')((n+1)Q+Q')} \]
(recall that $|QP' - PQ'| = 1$ as Farey neighbors.)   

Consequently,
\[ \frac{|\rho(t_{n+1}) - \rho(t_{n})|}{|\rho(t_{-\infty})-\rho_{\infty}|}
\frac{|t_{-\infty}-t_{\infty}|^{\alpha}}{|t_{n+1}-t_{n}|^{\alpha}} <
\frac {2(n+2)^{3\alpha-2}}{K_{1,P.\ref{prop:3,1}}^{\alpha}}  \; ,\]
where we used $Q<Q'<2Q$. By choosing $\alpha$ sufficiently close to $0$ we can 
get the ratio on the
right-side smaller than some number less than $1$, say $\beta$.
\end{proof}

Next, we want to generalize Lemma ~\ref{lem:n3,1} to $u,v$ arbitrary 
endpoints of the harmonic subdivision of $(P/Q, P'/Q')$. 

\begin{lem}\label{lem:n3,2}
Let $u,v$ be arbitrary two endpoints of the harmonic subdivision of 
$(P/Q,P'/Q')$. or $P/Q$ or $P'/Q'$. Then,
\[
|u-w||c(u)-c(w)|^{-\alpha}\leq K_{1}|P/Q-P'/Q'||c(P/Q)-c(P'/Q')|^{-\alpha}
\]
with $\alpha$ uniform and positive.
\end{lem}
\begin{proof}
First, we note that it suffices to prove the lemma when $u=u_{n}$ and 
$w=u_{n'}$ with $n\cdot n' \geq 0$. Indeed, in the situation when both
$u$ and $v$ are endpoints, but $n\cdot n' < 0$,
we can consider $n$ and $0$ as well as $n'$ and $0$ separately. Then, if we
sum up resulting inequalities and use convexity of $\frac{1}{x^{\alpha}}$, we
can infer the claim of the lemma.

If, for example, $u$ is $P/Q$, we can take the limit with $u_{n}$ where
$n$ tends to $+\infty$. By continuity of the 
rotation function, this would give us almost the estimate of the Lemma,
except that on the left-hand side $c(u)$ is replaced with $t_{0}$ which is
the upper endpoint of $\rho^{-1}(P/Q)$. However, this 
is stronger than the estimate claimed by the lemma.  

Furthermore, we can restrict our attention to $n,n' \geq 0$. We define $m$ by
$0<m=n'-n$. 

Again, we use Proposition ~\ref{prop:3,1}:
\[  \frac{|\rho(t_{n+m}) - \rho(t_{n})|}{|\rho(t_{-\infty})-\rho_{\infty}|}
\frac{|t_{-\infty}-t_{\infty}|^{\alpha}}{|t_{n+m}-t_{n}|^{\alpha}} \]
\[\leq\frac{n^{3\alpha} QQ'}
{K_{1,P.\ref{prop:3,1}}^{\alpha} (nQ+Q')((n+m)Q+Q')}
(\frac{1}{\sum_{k=n}^{n+m-1}1/k^{3}})^{\alpha} \leq\]
\[ \frac{2}{K_{1,P.\ref{prop:3,1}}^{\alpha}}\frac{m}{(m+n+1)(n+1)}
(\frac{1}{1/n^{2} - 1/(m+n)^{2}})^{\alpha} \leq\]
\[ \frac{2}{nK_{1,P.\ref{prop:3,1}}^{\alpha}}\frac{n^{2\alpha}(n+m)^{2\alpha}}
{m^{\alpha}(2n+m)^{\alpha}} \leq \frac{2^{2\alpha+1}}
{K_{1,P.\ref{prop:3,1}}^{\alpha}} n^{2\alpha-1} \; .\]

This expression is bounded by some $K_{1}$ if $\alpha\leq 1/2$. 
\end{proof}

\subsection{Global estimates}
Then, we let $u<w$ be arbitrary rational numbers from the unit interval. 
There is a unique simple path in the Farey tree from $u$ to $w$. It contains
the ``highest'' node $V$. This splits the path into two parts: from $u$ to
$V$ and from $V$ to $w$. 

Then, we define maps $\mu$ and $\nu$ on the Farey tree. Given a rational 
number $v$, $\mu(v)$ is the rational number that corresponds to the initial 
segment of the symbolic code of $v$ cut off at the last turning point 
(i.e. if the turning points are $ m_{1},\ldots, m_{k}$, the symbolic sequence 
of $\mu(v)$ is $a_{1},\ldots, a_{m_{k}}$.) 

Then, $\nu(v)$ is the mother of $\mu(v)$. 

Clearly, 
\[\mbox{deg}(\mu(v)) = \mbox{deg}(\nu(v)) = \mbox{deg}(v)-1\; .\]
Furthermore, $v$ lies strictly between $\mu(v)$ and $\nu(v)$, and $v$ is in 
fact an endpoint of the harmonic subdivision of $(\nu(v),\mu(v))$.
\footnote{In this section, if we write $(a,b)$, we perhaps mean $(b,a)$ when
$a>b$.}

\paragraph{The function $\zeta$.}
If $x$ and $y$ are two rationals from the unit interval, we define
\[ \zeta(x,y) := |x-y|\cdot |c(x)-c(y)|^{-\alpha}\; .\]

Lemma ~\ref{lem:n3,2} gives us a fundamental estimate
\begin{equation}\label{equ:n2,3}
\zeta(x,y) \leq K_{1} \zeta(\mu(x),\nu(x)) 
\end{equation}
provided $\mu(x) = \mu(y)$.

If we iterate $\mu$ and $\nu$,  we get nested sequence of growing
fundamental domains bounded by $\mu^{i}(x)$ and $\nu^{i}(x)$.
For any $i$ we have
\begin{equation}\label{equ:n2,4}
\zeta(\mu^{i}(x),\nu^{i}(x)) \leq \beta\zeta(\mu^{i+1}(x), \nu^{i+1}(x))
\end{equation}
by Lemma ~\ref{lem:n3,1}. 

Finally, we define a sequence $u_{0}=v$ and $u_{i}$ is equal to the 
greater of $\mu^{i}(v)$ and $\nu^{i}(v)$. Similarly, we define the 
sequence $w_{i}$ so that $w_{0}=w$ and $w_{i}$ is the minimum of 
$\mu_{i}(w)$ and $\nu_{i}(w)$.
Also, let $l$ be the largest so that $w_{l}\geq V$.

We want to bound $\zeta(u,u_{k})$. As $\zeta(v,v+v')$ is convex as a 
function of $v'$, and the sequence $u_{i}$ is growing, clearly
\begin{equation}\label{equ:n2,5}
 \zeta(u,u_{k}) \leq \sum_{i=1}^{k} \zeta(u_{i-1},u_{i})\; .
\end{equation}

Since $u_{i}$ and $u_{i+1}$ are contained between $\mu(u_{i+1})$
and $\nu(u_{i+1})$, we can estimate
\[ \zeta(u_{i},u_{i+1}) \leq K_{1}\zeta(\mu(u_{i+1}), \nu(u_{i+1})) \]
by inequality ~\ref{equ:n2,3}.
Then we iterate $\mu$ and $\nu$ on $\mu(u_{i+1})$ as many times as possible,
which at least $\max(k-i-1,0)$. By the repeated use of inequality ~\ref
{equ:n2,4}, we get
\[ \zeta(\mu(u_{i+1}),\nu(u_{i+1})) \leq \beta^{\max(0,k-i-1)}\zeta(0,1)\; .\]

Since $\zeta(0,1)=1$, we can finally get from estimate ~\ref{equ:n2,5}
\[ \zeta(u,u_{k}) \leq K_{2}\]
by simply adding up a geometric progression.

The same argument shows that 
\[ \zeta(w_{l},w) \leq K_{2} \; .\]   
Finally, by their definition $\mu(w_{l}) = \mu(u_{k})$ so that they
both are in the fundamental domain bounded by $\mu(u_{k})$ and $\nu(u_{k})$.

Again, we get 
\[ \zeta(u_{k}, w_{l}) \leq K_{1}\zeta(\mu(u_{k}),\nu(u_{k})) \]
from  estimate ~\ref{equ:n2,3},
and 
\[ \zeta(\mu(u_{k}),\nu(u_{k})) \leq 1\]
by the repeated use of estimate ~\ref{equ:n2,4}.

Since 
\[ \zeta(u,w) \leq \zeta(u,u_{k}) + \zeta(u_{k},w_{l}) + \zeta(w_{l},w) \]
by convexity, we have proved that  
\[\zeta(u,w) \leq K_{3}\]
with uniform $K_{3}$, which means the H\"{o}lder estimate. That is, we have
proven that
\[ |\rho(x) - \rho(y)|\leq K_{3} |x-y|^{\alpha} \]
provided that $x$ and $y$ are both centers of frequency-lockings.

\paragraph{The general Holder estimate.}
By continuity, we also get the same H\"{o}lder estimate if $x$ and $y$ 
belong to the closure of the set of centers, i.e. to the complement of the
union of interiors of all frequency-locking intervals. If $x<y$ are arbitrary,
we consider $x'$ which is the infimum of the set of centers which are between
$x$, and $y'$ is is the supremum of the same set. They are well-defined unless
$x$ and $y$ are in the same frequency-locking interval, in which case the 
estimate is evident. The H\"{o}lder estimate holds for $x'$ and $y'$. 
Moreover, $\rho(x')=\rho(x)$ and $\rho(y')=\rho(y)$, while $|x'-y'|\leq |x-y|$.
The H\"{o}lder estimate follows again.

\section{Hausdorff dimension}
We will prove Theorem B.
In this section, we will use the harmonic formalism of explained in 
Introduction section. That is, we have a one-to-one coding of fundamental
 domains 
in the parameter space by finite sequences of S-type symbols. The length of the
code will be called the {\em degree} of the corresponding fundamental domain.

The fundamental domain which corresponds to the code $(n_{1},\ldots,n_{r})$
will be denoted with ${\cal D}(n_{1},\ldots,n_{r})$.

Let $\Omega'$ be $T\setminus \bigcup_{w\in {\bf Q}} \rho^{-1}(w)$.

It is sufficient to prove that the Hausdorff dimension of 
\[ \Omega' \cap {\cal D}(n_{1},\ldots,n_{r}) \] 
satisfies our bounds for $r$ large enough. That means that we can assume that  
Proposition ~\ref{prop:1,3} holds and, consequently, use both claims of 
Proposition ~\ref{prop:3,1}.

\subsection{The estimate of HD($\Omega'$) from above}
 
Take a cover of $\Omega'$ which consists of all fundamental domains of degree
$r$.

Then 
\[ \sum_{n_{1},\ldots,n_{r}\in {\bf Z}\; ,\: |n_{r}|>k} 
|{\cal D}(n_{1},\ldots,n_{r})|^{\beta} \leq K_{2}(k) 
\sum_{n_{1},\ldots,n_{r-1}} |D(n_{1},\ldots,n_{r-1})|^{\beta} \]
if $\beta > 1/3$ as a consequence of scalings rules 
(Proposition ~\ref{prop:3,1}.) Here, $K_{2}(\cdot)$ is a uniform positive
function with limit $0$ at infinity.

Since for $|n_{r}|\leq k$ 
\[ K_{1} < \frac{|{\cal D}(n_{1},\ldots,n_{r})|}{|{\cal D}(n_{1},
\ldots,n_{r-1}|} \]
by Proposition ~\ref{prop:3,1}, by Young's inequality
\[ \sum_{n_{1},\ldots,n_{r}\in {\bf Z}\; ,\: |n_{r}|\leq k}
|{\cal D}(n_{1},\ldots,n_{r})|^{\beta} \leq k^{1-\beta}K_{1}^{\beta}
\sum_{n_{1},\ldots,n_{r-1}} |D(n_{1},\ldots,n_{r-1})|^{\beta}  \; .\]

Thus,
\[ \sum_{n_{1},\ldots,n_{r}\in {\bf Z}}
|{\cal D}(n_{1},\ldots,n_{r})|^{\beta} \leq (k^{1-\beta}K_{1}^{\beta}+K_{2}(k))
\sum_{n_{1},\ldots,n_{r-1}} |D(n_{1},\ldots,n_{r-1})|^{\beta} \; .\]

We claim that $k^{1-\beta}K_{1}^{\beta} + K_{2}(k)$ can be made less than 
$1$ by choosing $\beta$ sufficiently close to, but less than, one.
Indeed, remember that $K_{1}<1$. So, we first choose $k$ so large that 
$K_{2}(k)$ is less than $0.5(1-K_{1})$. Then, by adjusting $\beta$ we can make
$k^{1-\beta}K_{1}^{\beta}$ arbitrarily close to $K_{1}$. 

Since the diameters of the fundamental domains tend to $0$ with the degree, 
this $\beta$ is not larger than the Hausdorff dimension of $\Omega'$.

\subsection{Estimate from below by $1/3$.}
The proof is based on the following Frostman's Lemma which we borrowed from
\cite{felek}:
\begin{fact}\label{fa:n3,1}
Suppose that $\mu$ is a probabilistic Borel measure on the interval and that 
for $\mu$-a.e.  $x$
\[ \lim\inf_{\epsilon\rightarrow 0} \log(\mu(x-\epsilon, x+\epsilon))/
\log(\epsilon) \geq \lambda \; .\]
Then the Hausdorff dimension of $\mu$ is not less than $\lambda$.\footnote
{Frostman's Lemma remains true if the $\geq$ signs are replaced with $\leq$
signs. }
\end{fact}

Take $\eta < 1/3$. By the scalings rules it is clear that a number $k$ can be
found independently of the fundamental domain ${\cal D}(n_{1},\ldots,n_{r})$
so that
\[ \sum_{n_{1},\ldots,n_{r}\in {\bf Z}\; ,\: |n_{r}|\leq k} |{\cal D}
(n_{1}, \ldots, n_{r})|^{\eta} \geq |{\cal D}(n_{1},\ldots,n_{r-1})|^{\eta}\; 
.\]

We now define $\mu$ as a limit of probabilistic measures. The measure
$\mu_{0}$ is just the Lebesgue measure on $\rho^{-1}(0,1)$ properly scaled. 
To obtain $\mu_{i+1}$, we consider all fundamental
domains of degree $i$. If ${\cal D}(n_{1},\ldots,n_{i})$ is one of those,
the density of $\mu_{i+1}$ with respect to $\mu_{i}$ on 
 ${\cal D}(n_{1},\ldots,n_{i})$ equals
\[ \frac{|{\cal D}(n_{1},\ldots,n_{i+1}|^{\eta-1}}{\sum_{n_{i+1}=-k}^
{n_{i+1}=k}|{\cal D}(n_{1},\ldots,n_{i+1})|^{\eta}} \]
on fundamental domains of the harmonic subdivision with $|n_{j+1}|\leq k$,
and is zero the mode-locking intervals which belong to the endpoints of the 
subdivision. The sequence has a limit which is supported on a set $\Omega_{0}$ 
contained in $\Omega'$. The set $\Omega_{0}$ consists of the preimages in the 
parameter space
of all irrationals with harmonic code symbols $S(n,n+1)$ with $n,n+1$ bounded 
by $k$ as to absolute value.
Moreover, by our choice of $\eta$ 
we see that the denominator is not less that 
$|{\cal D}(n_{1},\ldots,n_{i-1})|^{\eta}$, thus the the density is not greater
than 
\[ \frac{|{\cal D}(n_{1},\ldots,n_{i+1})|^{\eta-1}}{|{\cal D}(n_{1},\ldots,
n_{i})|^{\eta}}\; , \]
therefore 
\[ \frac{\mu_{i+1}({\cal D}(n_{1},\ldots,n_{i+1}))}{\mu_{i}
({\cal D}(n_{1},\ldots,n_{i}))} \leq 
\frac{|{\cal D}(n_{1},\ldots,n_{i+1})|^{\eta}}
{|{\cal D}(n_{1},\ldots,n_{i})|^{\eta}} \; .\]

So, by induction,
\[ \mu_{i+1}({\cal D}(n_{1},\ldots,n_{i+1})) \leq 
|{\cal D}(n_{1},\ldots,n_{i+1})|^{\eta}\; , \]
and since clearly 
\[ \mu({\cal D}(n_{1},\ldots,n_{i+1})) = 
\mu_{i+1}({\cal D}(n_{1},\ldots,n_{i+1}))\; ,\]
the same estimate holds for $\mu$ itself. 

Take any small $\epsilon>0$ and an $x\in \Omega_{0}$ and look for the largest 
$r$ so that
\[ \mu((x-\epsilon,x+\epsilon) \setminus {\cal D}(n_{1},\ldots,n_{r}))=0 \]
for some $n_{1},\ldots,n_{r}$. Note that a finite $r$ with this property always
exists by topology.

Then 
\[ \mu((x-\epsilon,x+\epsilon))\leq |{\cal D}(n_{1},\ldots,n_{r})|^{\eta} \]
and 
$2\epsilon$ is greater than the length of the some gap between domains of
the harmonic subdivision of ${\cal D}(n_{1},\ldots,n_{r})$ which have non-
zero measure. Indeed, by the definition of $r$, the interval $(x-\epsilon,
x+\epsilon)$ must be straddled between at least two such domains. If the
size of the gap is denoted with $\gamma$, we get
\[ \frac{\log(\mu(x-\epsilon,x+\epsilon))}{\log(\epsilon)}\geq
\eta\frac{\log(|{\cal D}(n_{1},\ldots,n_{r})|)}{\log(\gamma)} \; .\]

By Proposition ~\ref{prop:3,1}, the lengths of both domains
and $|{\cal D}(n_{1},\ldots,n_{r})|$ are all related by uniform constants.
But $\gamma$ is not much smaller than either of them as a result of the
estimates of ~\cite{rat}. So the 
logarithms differ by a bounded amount and their ratio tends to $1$ as
$\epsilon$ shrinks to $0$. 

So, if we pass to the limit with $\epsilon\rightarrow 0$ we see that the 
assumptions of Frostman's Lemma are satisfied, therefore the Hausdorff 
dimension of $\Omega_{0}$, which is larger than $\Omega'$, is at least
$\eta$. But $\eta$ could have be chosen anything less than $1/3$. So, the
estimate follows.  

\subsection{Justification of Remark B}
Remark B follows from the fact that our estimates can be made uniform 
with respect to the family if we consider Farey domains of sufficiently 
large degree. Namely, it can be checked that all
our estimates ultimately depend on bounded geometry and the Distortion
Lemma. 

We claim that as far as bounded geometry is concerned, for every family
a uniform $K_{2}$ can be chosen so that if $k$ which occurs in the statement
of bounded geometry is greater than $K_{2}$, the estimate becomes uniform with
respect to the family. This is a known fact which was observed and explained
in \cite{miszczu}. The reason is that ``cross-ratio inequalities'' 
(see \cite{rat}) which are the source of bounded geometry estimates become 
uniform if applied to very short intervals. This fact can be seen immediately
from the ``pure singularity property'' of \cite{poincare}, but also follows
easily from the much simpler Corollary to Proposition ~\ref{prop:1,3}.     

The Distortion Lemma also becomes uniform if applied to maps with rotation
number of degree sufficiently large. This follows from the Corollary to 
Proposition ~\ref{prop:1,3}. Indeed, the nonlinearity of the $N$ part is 
uniformly bounded by the classical estimate of \cite{guke}, while the 
other parts are also bounded if $J$ is small, simply because the joint length
of $i_{1}$ or $i_{2}$ images of $J$ is small.

Thus, for every family the scalings rules become uniform with respect to the
family on Farey domains of large degree. We want to emphasize that how large
that degree should be depends on the family. However, since the Hausdorff 
dimension is an asymptotic quantity which we bound from scalings rules, it 
will be uniformly bounded away from $1$.

\end{document}